\documentclass[12pt]{article}
\usepackage[fleqn]{amsmath}
\usepackage[cp1251]{inputenc}
\usepackage[english,russian]{babel}
\usepackage{amsmath,amsfonts,amssymb,graphicx,makeidx,amsthm}
\newtheorem{teo}{Теорема}

\begin{document}

\begin{center}
 \textbf{Э.Л. Шишкина}
\end{center}

\begin{center}
 \textbf{ВЕСОВЫЕ СФЕРИЧЕСКИЕ СРЕДНИЕ, ПОРОЖДЕННЫЕ ОБОБЩЕННЫМ СДВИГОМ, И ОБЩЕЕ УРАВНЕНИЕ ЭЙЛЕРА-ПУАССОНА-ДАРБУ.\\ ЧАСТЬ 1.}
\end{center}

\renewcommand{\contentsname}{Содержание}

\tableofcontents

\newpage

\section{Определение весового сферического среднего}

В этой статье мы приведем результаты, полученные для  весовых сферических средних, порожденных обобщенным сдвигом, и решение задачи Коши для   общего уравнения Эйлера-Пуассона-Дарбу.

Начнём с определения классического сферического среднего. Рассмотрим евклидово пространство $\mathbb{R}^n$. Пусть $x=(x_{1} ,\ldots ,x_{n} )\in\mathbb{R}^n$ и $|x|=\sum\limits_{j=1}^{n} x_{j}^{2} $.
Классическое сферическое среднее  имеет вид
\begin{equation}\label{ShSr}
M(x,r,u)=\frac{1}{|S_{n}(1)|} \int\limits_{S_{n}(1)} u(x+\beta r)dS,
\end{equation}
где $S_{n}(1) =\left\{\; x\in\mathbb{R}^n:|x|=1\; \right\}$ --- сфера единичного радиуса с центром в начале координат, $r\ge 0$, $f(x+\beta r)=f(x_{1} +\beta _{1} r,\ldots ,x_{n} +\beta _{n} r)$, $|S_{n}(1)|=\int\limits_{S_{n} } dS=2\frac{\pi ^{\frac{n}{2} } }{\Gamma \left(\frac{n}{2} \right)} $ --- площадь сферы $S_{n}(1)$, $\Gamma(\alpha)$ --- гамма функция Эйлера, $dS$ --- элемент поверхности $S_{n}(1)$.

 Операторы \eqref{ShSr} естественным образом возникают в математике при решении различных задач для уравнений гиперболического типа. Классический метод сведения волнового уравнения 
 \begin{equation}\label{Voln}
 \sum\limits^n_{i=1}\frac{{\partial }^2u}{\partial x^2_i}=\frac{{\partial }^2u}{\partial t^2},\qquad u=u(x_1,...,x_n,t)
\end{equation}
в случае $n\geq 2$ пространственных переменных  к волновому уравнению с одной пространственной переменной посредством оператора \eqref{ShSr} при $n=3$ изложен С. Д. Пуассоном  (Sim${\rm \acute{e}}$on Denis Poisson) в \cite{Poisson0} (см., также, \cite{Adamar} стр. 50, \cite{Bers} стр. 22). В случае произвольного $n\geq 2$ решение \eqref{Voln} выражено через \eqref{ShSr} в \cite{john} на стр. 36. Теорема о среднем  \eqref{ShSr} для регулярной гармонической функции приведена в
  \cite{kurantgil2} на стр. 278. Л. Асгейрссон (Leifur $\acute{{\rm A}}$sgeirsson) в \cite{Asgeirsson} доказал теорему о сферических средних для ультрагиперболического уравнения вида
\begin{equation}\label{UG1}
\sum\limits_{j=1}^{n}\frac{\partial^2 u}{\partial x_j^2}=\sum\limits_{j=1}^{n}\frac{\partial^2}{\partial y_j^2u},\qquad u=u(x_1,...,x_n,y_1,...,y_n).
\end{equation}
 Р. Курант (Richard Courant) в книге \cite{kurantgil2}, стр. 473 (см., также, \cite{kurant2}, стр. 738) применяет теорему Асгейрссона к решению задачи Коши для волнового уравнения. Ф. Йон (Fritz John) \cite{john}, стр. 87 доказывает единственность решения задачи О.Л. Коши (Augustin Louis Cauchy) вида
$$
\sum\limits^n_{i=1}\frac{{\partial }^2u}{\partial x^2_i}=\frac{{\partial }^2u}{\partial t^2}+\frac{k}{t}\frac{\partial u}{\partial t},\quad u=u(x,t),\quad x\in\mathbb{R}^n,\quad t>0,\quad 0<k<\infty,
$$
$$
u(x,0)=f(x),\qquad \frac{\partial u}{\partial t}\biggr|_{t=0}=0
$$
  и получает решение этой задачи. Л. Хёрмандер (Lars H\"ormander) в \cite{xer}, стр. 222 доказал уточнение теоремы Асгейрссона, при помощи которого затем, получена формула Г. Кирхгофа (Gustav Kirchhoff). Исследования вопросов корректности постановки и разрешимости задач для ультрагиперболического уравнения \eqref{UG1} при помощи теоремы Асгейрссона приведены в \cite{Blago1}--\cite{Radzikowski}.

Большой интерес у различных исследователей вызывает обобщение сферического среднего \eqref{ShSr}. Так, в работе \cite{Weinstein2} рассмотрено сферическое среднее в пространстве с отрицательной кривизной, в \cite{Dunkl} и \cite{Elouadih} рассмотрено обобщение сферического среднего, порожденное оператором преобразования С. Дункла (Charles Dunkl). В данной работе будет рассмотрено обобщение сферического среднего на случай, когда вместо ультрагиперболического уравнения \eqref{UG1} рассматривается уравнение типа  ультрагиперболического, но в котором вместо второй производной по каждому аргументу применяется дифференциальный оператор Бесселя вида
$$
(B_\nu)_x=\frac{\partial^2}{\partial x^2}+\frac{\nu}{x}\frac{\partial}{\partial x}.
$$
Такое уравнение называется В-ультрагиперболическим и имеет вид
\begin{equation}\label{UG22}
\sum\limits_{j=1}^{n}(B_\nu)_{x_j}u=\sum\limits_{j=1}^{n}(B_\nu)_{y_j}u,\qquad u=u(x_1,...,x_n,y_1,...,y_n).
\end{equation}
При построении  сферического среднего для уравнения \eqref{UG22}  вместо обычного сдвига применяется многомерный обобщенный сдвиг. О схеме конструирования одномерного обобщенного сдвига, предложенной Ж. Ф. Дельсартом (Jean Fr${\rm \acute{e}}$d${\rm \acute{e}}$ric Delsarte) кратко приведены сведения в пункте 1.1 этой статьи. В этом же пункте  приведен обобщённый сдвиг, порожденный оператором Бесселя, при помощи которого и конструируется рассматриваемое нами в статье весовое сферическое среднее. Некоторые формулы для интегралов по частям сфер приведены в  пункте 1.2. В пункте 1.3 приведены определения многомерного обобщенного сдвига, весового сферического, итерированного весового сферического и некоторые формулы для этих операторов.

\subsection{Операторы обобщённого сдвига с точки зрения Дельсарта}

Рассмотрим путь обобщения оператора сдвига, который был предложен Ж. Дельсартом (Jean Delsarte) в \cite{Delsarte1}--\cite{Delsarte4} (см. также \cite{Marchenko}--\cite{McGregor}).

Если $f$ --- функция, определенная на вещественной оси, то оператор сдвига $T^y$, $y\in\mathbb{R}$, определяется равенством
\begin{equation}\label{Shift}
T^y_xf(x)=f(x+y).
\end{equation}
Пусть теперь $f\in C^\infty(\mathbb{R})$. Подход Ж. Дельсарта заключался в нахождении обобщения формулы Тэйлора
\begin{equation}\label{EQ00}
T^y_x f(x)=f(x+y)=\sum\limits_{n=0}^\infty \frac{y^n}{n!}\left(\frac{d}{dx}\right)^n\, f(x),
\end{equation}
которая даёт разложение оператора сдвига $T^y$ по степеням
оператора дифференцирования $\frac{d}{dx}$. Для сдвига \eqref{Shift}, опираясь на \eqref{EQ00} Дельсарт сопоставлял функцию $\varphi_n(y)=\frac{y^n}{n!}$ дифференциальному оператору $L_x=\frac{d}{dx}$ в некотором специальном смысле. А именно, он исходил из того, что  решение $\varphi(y,\lambda)$, $y\in\mathbb{R}$, $\lambda\in\mathbb{C}$ задачи
\begin{equation}\label{EQ01}
L_y\varphi=\lambda\varphi,\qquad \varphi(0,\lambda)=1
\end{equation}
есть $\varphi(y,\lambda)=e^{\lambda y}$. Для каждого вещественного $y$ эта функция является целой функцией $\lambda$ и
\begin{equation}\label{EQ02}
\varphi(y,\lambda)=\sum\limits_{n=0}^\infty \varphi_n(y)\lambda^n\qquad {\text{или}}\qquad e^{\lambda y}=\sum\limits_{n=0}^\infty \frac{y^n}{n!}\,\lambda^n.
\end{equation}
Функции $\varphi_n(y)=\frac{y^n}{n!}$, $n=0,1,2,...$ удовлетворяют условиям
$$
L_y \varphi_0=0,\qquad \varphi_0(1)=1,
$$
$$
L_y \varphi_n=\varphi_{n-1},\qquad \varphi_n(0)=0,\qquad n=1,2,...
$$
C учетом \eqref{EQ01}, Дельcарт обобщает формулу Тэйлора \eqref{EQ00}  следующим образом:
$$
T^y_x f(x)=\sum\limits_{n=0}^\infty \varphi_n(y)(L_x)^nf(x),
$$
где $L_x$ --- некоторый оператор. Очевидно, что поскольку $L_y \varphi_0=0$ и $L_y \varphi_n=\varphi_{n-1}$, то формально
$$
L_yT^y_x f(x)=\sum\limits_{n=0}^\infty L_y\varphi_n(y)(L_x)^nf(x)=\sum\limits_{n=1}^\infty \varphi_{n-1}(y)(L_x)^nf(x)=$$
$$=\sum\limits_{n=0}^\infty \varphi_n(y)(L_x)^{n+1}f(x)=L_xT^y f(x),
$$
то есть  $T^y_x f(x)$ формально удовлетворяет дифференциальному
уравнению
\begin{equation}\label{Zad1}
L_xT^y_x f(x)=L_yT^y_x f(x)
\end{equation}
при начальных условиях
\begin{equation}\label{Zad2}
T^y_xf(x)|_{y=0}=f(x),\qquad \frac{\partial}{\partial y}T^y_xf(x)\biggr|_{y=0}=0.
\end{equation}
Операторы $T^y_x$ Дельсарт назвал операторами обобщённого сдвига и установил для них ряд свойств.
Обобщенный сдвиг является одним из примеров оператора преобразования (см. \cite{S42}--\cite{S5}).

В  работе \cite{Delsarte2} Дельсарт  подробно изучает оператор обобщённого
сдвига в том случае, когда
$$
L_x=(B_\gamma)_x=\frac{\partial^2}{\partial x^2}+\frac{\gamma}{x}\frac{\partial}{\partial x}
$$
есть дифференциальный оператор Бесселя. Такой обобщенный сдвиг также подробно рассматривается в статье \cite{levitan1}. %Краткие сведения об этом обобщенном сдвиге приведены в следующем пункте.

%\subsection{Обобщённый сдвиг, порожденный оператором Бесселя}

Рассмотрим обобщенный сдвиг, порожденный оператором Бесселя $B_\gamma$ и приведем некоторые его свойства. Пусть $x\geq0$, $y\geq0$. Подставляя в \eqref{Zad1}--\eqref{Zad2} оператор
\begin{equation}\label{Bess}
(B_\gamma)_x=\frac{\partial^2}{\partial x^2}+\frac{\gamma}{x}\frac{\partial}{\partial x},\qquad \gamma>0,
\end{equation}
 получим  задачу Коши вида
\begin{equation}\label{Zad3}
(B_\gamma)_xT^y_x f(x)=(B_\gamma)_yT^y_xf(x),
\end{equation}
\begin{equation}\label{Zad4}
T^y_xf(x)|_{y=0}=f(x),\qquad \frac{\partial}{\partial y}T^y_xf(x)\biggr|_{y=0}=0.
\end{equation}
Если $f\in C^2(0,\infty)$, то единственное решение задачи \eqref{Zad3}--\eqref{Zad4} дается формулой (см. \cite{Delsarte2}, \cite{levitan1})
\begin{equation}\label{Sdvog0}
T^y_xf(x)=C(\gamma)\int\limits_0^\pi
f(\sqrt{x^2+y^2-2xy\cos{\varphi}})\sin^{\gamma-1}{\varphi}d\varphi,
\end{equation}
$$
C(\gamma)=\left(\int\limits_0^\pi\sin^{\gamma-1}{\varphi}d\varphi\right)^{-1}=\frac{\Gamma\left(\frac{\gamma+1}{2}\right)}{\sqrt{\pi}\,\,\Gamma\left(\frac{\gamma}{2}\right)}.
$$
Отметим, что поскольку  $T^y_xf$ зависит не только от $x$ и $y$, но еще и от $\gamma$, то наряду с обозначением $T^y_xf$ далее будем использовать также обозначение $\,^\gamma T^y_xf$.
В силу второго условия в \eqref{Zad4} будем рассматривать обобщенный сдвиг $T^y_xf$  только таких функций $f$, для которых $\frac{\partial}{\partial y}f(y)\biggr|_{y=0}=0$. Исходя из этого, продолжение на значения $y<0$ рассматриваемых функций будем осуществлять четным образом.

Обобщенный сдвиг можно записать в виде
\begin{equation}\label{Sdvig}
T^y_xf(x)=\frac{2^\gamma C(\gamma)}{(4xy)^{\gamma-1}}\int\limits_{|x-y|}^{x+y}zf(z)[(z^2-(x-y)^2)((x+y)^2-z^2)]^{\frac{\gamma}{2}-1}dz.
\end{equation}
Используя \eqref{Sdvig}, легко доказать свойство самосопряженности обобщенного сдвига, доказанное в \cite{levitan1} другим способом:
\begin{equation}\label{Samos}
\int\limits_0^{+\infty}T^y_xf(x)g(y)y^\gamma dy
=\int\limits_0^{+\infty}f(y)T^y_xg(x)y^\gamma dy.
\end{equation}
Как уже отмечалось, мы рассматриваем действие обобщенного сдвига только на четные функции, поэтому в формуле \eqref{Samos} интегрирование ведется по полуоси, а не по всей оси. Кроме того, в формуле \eqref{Samos} необходимо присутствие веса $y^\gamma$. Необходимость наличие степенного веса $y^\gamma$ можно объяснить следующим образом.
Пусть $x=(x_1,...,x_n)\in\mathbb{R}^n$. Рассмотрим оператор Лапласа, действующий по $n$ переменным:
$$ \Delta  = \sum_{i=1}^n \frac {\partial^2 }{\partial x^2_i}.
$$
В сферических координатах $x=r\theta$, $r\geq 0$, $\theta\in\mathbb{R}_n$ оператор Лапласа примет вид
$$
\Delta ={\frac {\partial ^{2}}{\partial r^{2}}}+{\frac {n-1}{r}}{\frac {\partial }{\partial r}}+{\frac {1}{r^{2}}}\Delta _{S^{n-1}},
$$
где $\Delta_{S^{n-1}}$ --- оператор Лапласа-Бельтрами, действующий по правилу
$$
\Delta_{S^{n-1}}f(x)=\Delta f\left(\frac{x}{|x|}\right).
$$
Если функция $n$ переменных зависит только от $|x|{=}\sqrt{x_1^2+...+x_n^2}{=}r$, то, очевидно, что применяя к ней оператора Лапласа в сферических координатах мы получаем оператор Бесселя \eqref{Bess} при натуральном $\gamma=n-1$:
$$
\Delta f(r)=\left({\frac {\partial ^{2}}{\partial r^{2}}}+{\frac {n-1}{r}}{\frac {\partial }{\partial r}}\right)f(r)=(B_{n-1})_rf(r).
$$
Таким образом, в случае сферически симметричной функции $f=f(|x|)=f(r)$ оператор Бесселя при натуральном $\gamma=n-1$ получается при переходе к сферическим координатам в операторе Лапласа, при этом $\gamma=n-1$ --- размерность пространства $\mathbb{R}^n$, в котором рассматривалась функция $f=f(|x|)=f(r)$. Рассмотрим интеграл по $\mathbb{R}^n$ от функции $f=f(|x|)=f(r)$ и произведем в рассматриваемом интеграле переход к сферическими координатам $x=r \theta$, получим
$$
\int\limits_{\mathbb{R}^n}f(|x|)dx=\int\limits_0^{+\infty} f(r)r^{n-1}dr\int\limits_{S_1(n)}dS,
$$
где $S_r(n)$ --- сфера радиуса $1$  с центром  начале координат в $\mathbb{R}^n$. Поскольку $\int\limits_{S_1(n)}dS{=}\frac{2\pi^{\frac{n}{2}}}{\Gamma\left(\frac{n}{2}\right)}$ (см. \cite{Coscia}), то
$$
\int\limits_{\mathbb{R}^n}f(|x|)dx=\frac{2\pi^{\frac{n}{2}}}{\Gamma\left(\frac{n}{2}\right)}\int\limits_0^{+\infty} f(r)r^{n-1}dr.
$$
Таким образом, когда  мы рассматриваем интеграл по $x$ от обобщенного сдвига $\,^{n-1}T^y_x$, соответствующего оператору $B_{n-1}$, то нужно умножать подынтегральную функцию на степенной вес $x^{n-1}$. Если рассматривать обобщенной сдвиг $\,^{\gamma}T^y_x$ для произвольного положительного вещественного число $\gamma$, то беря интеграл по $x$ , содержащий такой обобщенный сдвиг нужно умножать подынтегральную функцию на степенной вес $x^{\gamma}$.

%\vskip 1cm

Помимо  \eqref{Zad4} для обобщенного сдвига справедливы следующие свойства (см. \cite{levitan1}).
\begin{enumerate}
  \item Линейность и однородность: $$T^y_x[af(x)+bg(x)]=aT^y_xf(x)+bT^yg(x),\qquad a,b\in\mathbb{R}.$$
  \item Положительность: $T^y_xf(x)\geq 0$, если $f(x)\geq 0.$
  \item $T^y_x[1]=1.$
  \item Если $f(x)\equiv 0$ для $x\geq a$, то $T^y_x F(x)\equiv 0$ для $|x-y|\geq a$.
  \item Если последовательность непрерывных функций $f_n(x)$ сходится равномерно в каждом конечном интервале к $f(x)$, то последовательность функций от двух переменных $T^y_x f_n(x)$  сходится равномерно в каждой конечной области к функции $T^y_x f(x)$.
       \item Оператор $T^y_x$ ограничен: $$|T^y_xf(x)|\leq T^y_x|f(x)|\leq \sup\limits_{x\geq 0}|f(x)|.$$
       \item  Переместительность операторов $T^y_x$: $$T^y_xT^z_x f(x)=T^z_xT^y_xf(x).$$
        \item Ассоциативность операторов $T^y_x$: $$T^z_yT^y_xf(x)=T^z_xT^y_xf(x).$$
\end{enumerate}

Рассмотрим функцию $j_\nu$, определяемую формулой (см. \cite{kipr}, стр. 10)
\begin{equation}\label{FBess}
j_\nu(x) ={2^\nu\Gamma(\nu+1)\over x^\nu}\,\,J_\nu(x),
\end{equation}
где  $J_\nu$ --- функция Бесселя первого рода.
Для \eqref{FBess} справедливо равенство (см., например, \cite{levitan1})
\begin{equation}\label{RavDBes}
 T^y_x j_{\frac{\gamma-1}{2}}(x)=j_{\frac{\gamma-1}{2}}(x)\,j_{\frac{\gamma-1}{2}}(y).
\end{equation}

Приведем еще одну полезную формулу, представляющую обобщенный сдвиг от степенной функции $x^\alpha$:
$$
T^y x^\alpha =|x-y|^\alpha
\,_2F_1\left(-\frac{\alpha}{2},\frac{\gamma}{2},\gamma;-\frac{4xy}{(x-y)^2}\right),
$$
где $\,_2F_1$ --- гипергеометрическая функция Гаусса (см., например, \cite{Abramowitz}, стр. 370).

Очевидно, что при $\gamma=0$ в \eqref{Bess} мы получаем вторую производную. Поэтому при $\gamma=0$ задача Коши \eqref{Zad3}--\eqref{Zad4} примет вид
\begin{equation}\label{Zad5}
\frac{\partial^2}{\partial x^2}T^y_x f(x)=\frac{\partial^2}{\partial y^2}T^y_xf(x),
\end{equation}
\begin{equation}\label{Zad6}
T^y_xf(x)|_{y=0}=f(x),\qquad \frac{\partial}{\partial y}T^y_xf(x)\biggr|_{y=0}=0
\end{equation}
и $T^y_x$  будет представлять собой
\begin{equation}\label{RegSdv}
T^y_xf(x)=\frac{f(x+y)-f(x-y)}{2}.%=\frac{1}{2}\int\limits^{x+y}_{x-y}f'(z)dz
\end{equation}
Логично ожидать, что для четной функции $f(x)$ при $\gamma=0$ формула \eqref{Sdvog0} даст \eqref{RegSdv}, однако нельзя подставить $\gamma=0$ в интеграл
$$
\int\limits_0^\pi
f(\sqrt{x^2+y^2-2xy\cos{\varphi}})\sin^{\gamma-1}{\varphi}d\varphi,
$$
так как он будет расходящимся. Поэтому рассмотрим интеграл по $y$ от $0$ до $+\infty$ от $^\gamma T^y_xf(x)$, $\gamma>0$ с весом $y^\gamma$, применим формулу  \eqref{Samos} и свойство 3, а затем, предполагая, что $f$ финитна, следовательно, интеграл $\int\limits_0^{+\infty}f(y)y^\gamma dy$ сходится равномерно при $\gamma>0$, перейдем к пределу при $\gamma\rightarrow 0$:
$$
\int\limits_0^{+\infty}T^y_xf(x) y^\gamma dy= \int\limits_0^{+\infty}f(y)T^y_x[1] y^\gamma dy=\int\limits_0^{+\infty}f(y)y^\gamma dy,
$$
$$
\lim\limits_{\gamma\rightarrow 0}\int\limits_0^{+\infty}T^y_xf(x) y^\gamma dy=\lim\limits_{\gamma\rightarrow 0}\int\limits_0^{+\infty}f(y)y^\gamma dy=\int\limits_0^{+\infty}f(y)dy.
$$
С другой стороны, если $f$ --- четная финитная функция, то
$$
\lim\limits_{\gamma\rightarrow 0}\int\limits_0^{+\infty}\frac{f(x+y)-f(x-y)}{2}y^\gamma dy=\frac{1}{2}\left(\int\limits_0^{+\infty}f(x+y)dy-\int\limits_0^{+\infty}f(x-y)dy\right)=
$$
$$
=\frac{1}{2}\left(\int\limits_y^{+\infty}f(t)dt-\int\limits_{-y}^{+\infty}f(t)dt\right)=\frac{1}{2}\left(\int\limits_y^{+\infty}f(t)dt+\int\limits_{-\infty}^{y}f(t)dt\right)=
$$
$$
=\frac{1}{2}\int\limits_{-\infty}^{+\infty}f(t)dt=\int\limits_0^{+\infty}f(t)dt,
$$
следовательно, для  четной финитной функции $f$ справедливо равенство
$$
\lim\limits_{\gamma\rightarrow 0}\int\limits_0^{+\infty}T^y_xf(x) y^\gamma dy=\lim\limits_{\gamma\rightarrow 0}\int\limits_0^{+\infty}\frac{f(x+y)-f(x-y)}{2}y^\gamma dy=\int\limits_0^{+\infty}f(t)dt.
$$
Если в определениях и утверждениях классического гармонического анализа заменить обычный сдвиг обобщенным, то получится, так называемый, \textbf{весовой гармонический анализ} (см. \cite{kipr}, \cite{LyahovKniga1}--\cite{platonov1})  из-за наличия степенного веса $x^\gamma$ под знаком интегралов (см., например, формулу \eqref{Samos}).

\subsection{Часть евклидова пространства и функции, определенные в ней. Интеграл по части сферы и интеграл по части шара}

Из предыдущего пункта следует, что для того чтобы действие обобщенного сдвига было определено корректно, во-первых, целесообразно рассматривать положительные значения переменных, чтобы присутствующий вес не являлся бы многозначной функцией, во-вторых, для рассматриваемых дифференцируемых функций должно быть выполнено условие равенства нулю первой производной в начале координат.

Исходя из вышесказанного, рассмотрим часть  пространства $\mathbb{R}^n$ вида
$$
\mathbb{R}^n_+{=}\{x{=}(x_1,\ldots,x_n)\in\mathbb{R}^n,\,\,\, x_1{>}0,\ldots, x_n{>}0\}
$$
и $\Omega$ открытое множество в $\mathbb{R}^n$, симметричное относительно каждой гиперплоскости $x_i{=}0$, $i=1,...,n$, $\Omega_+=\Omega\cap{\mathbb{R}}\,^n_+$ и $\overline{\Omega}_+=\Omega\cap\overline{\mathbb{R}}\,^n_+$ где $$
\overline{\mathbb{R}}\,^n_+{=}\{x{=}(x_1,\ldots,x_n)\in\mathbb{R}^n,\,\,\, x_1{\geq}0,\ldots, x_n{\geq}0\}.
$$ Имеем $\Omega_+\subseteq{\mathbb{R}}\,^n_+$ и $\overline{\Omega}_+\subseteq\overline{\mathbb{R}}\,^n_+$.
Мы рассмотрим множество $C^m(\Omega_+)$, состоящее из $m$ раз дифференцируемых на on $\Omega_+$ функций.
Через $C^m(\overline{\Omega}_+)$ обозначим подмножество функций из $C^m(\Omega_+)$ таких, что все производные этих функций по $x_i$ для любого $i=1,...,n$ непрерывно продолжаются на $x_i{=}0$. Функции $f\in C^m(\overline{\Omega}_+)$ мы будем называть \textbf{четными по переменным}  $x_i$, $i=1,...,n$ если $\frac{\partial^{2k+1}f}{\partial x_i^{2k+1}}\biggr|_{x=0}=0$ для всех неотрицательных целых $k\leq m$ (см. \cite{kipr}, стр. 21). Класс $C^m_{ev}(\overline{\Omega}_+)$ состоит из функций из  $C^m(\overline{\Omega}_+)$, четных по каждой из своих переменных $x_i$, $i=1,...,n$. Тогда
$$
C^\infty_{ev}(\overline{\Omega}_+)=\bigcap\limits_{m=0}^{\infty}C^m_{ev}(\overline{\Omega}_+).
$$

Переходя к многомерному случаю нам необходимо также перейти от числа $\gamma$ в \eqref{Bess} к вектору $(\gamma_1,...,\gamma_n)$. А именно, будем рассматривать мультииндекс $\gamma{=}(\gamma_1,{\ldots},\gamma_{n})$, состоящий из положительных фиксированных чисел $\gamma_i>0$, $i{=}1,{...},n$ и   положим $|\gamma|{=}\gamma_1{+}{\ldots}{+}\gamma_{n}.$

Часть шара  $|x|\leq r$, $|x|=\sqrt{x_1^2+...+x_n^2}$, принадлежащую $\mathbb{R}^n_+$, будем обозначать
$B_r^+(n)$. Граница $B_r^+(n)$ состоит из части сферы $S_r^+(n){=}\{x\in\mathbb{R}^n_+:|x|{=}r\}$ и из частей координатных гиперплоскостей $x_i{=}0$, $i{=}1,{\ldots},n$, таких что $|x^i|\leq r$\,.
Для определения весового сферического среднего нам потребуется рассматривать интеграл по  $S_r^+(n)$ с весом
$$x^\gamma=\prod\limits_{i=1}^nx_i^{\gamma_i}$$ вида
$$
\int\limits_{S_r^+(n)}u(x)x^\gamma dS_r,
$$
где $dS_r$ --- элемент поверхности $S_r^+(n)$. Легко видеть, что справедливы формулы
\begin{equation}\label{ArBall}
\int\limits_{B_r^+(n)}u(x)x^\gamma dx=r^{n+|\gamma|}\int\limits_{B_1^+(n)}u(rx)x^\gamma dx,
\end{equation}
и
\begin{equation}\label{ArSph}
\int\limits_{S_r^+(n)}u(x)x^\gamma dS_r=r^{n+|\gamma|-1}\int\limits_{S_1^+(n)}u(rx)x^\gamma dS,
\end{equation}
где $dS$ --- элемент поверхности $S_1^+(n)$

Для  интегрируемой по $B_r^+(n) $ с весом $x^\gamma$ функции $f(x)$ и для непрерывной при $t\in[0,\infty)$  функции $g(t)$,  имеют место  формулы (см., например, \cite{LPSh1})
 \begin{equation}\label{Ball}
\int\limits_{B_r^+(n)}g(|x|)f(x)\,x^\gamma\, dx=\int\limits_0^r g(\lambda)\lambda^{n+|\gamma|-1}d\lambda \int\limits_{S_1^+(n)}f(\lambda x)x^\gamma d\omega\,,
 \end{equation}
 и при $g(x)=1$ и непрерывной на $B_r^+(n)$ функции $f(z)$
\begin{equation}\label{Ball1}
\int\limits_{S^+_1(n)}f(r x)x^\gamma dS=r^{1-n-|\gamma|}\frac{d}{dr}\int\limits_{B_r^+(n)}f(z)z^\gamma dz.
 \end{equation}
 Кроме того, нам знать значение интеграла $\int\limits_{S^+_1(n)} x^{\gamma}dS$, представляющего площадь взвешенной части сферы $S_1^+(n)$. Для нахождения этого интеграла используем формула, которая доказывается методом математической индукции (см. \cite{Coscia}):
\begin{equation}\label{GammaS}
    \mathop{\int...\int}_{y_1+...+y_n\leq 1,y_1\geq 0,...,y_n\geq 0}y^{\alpha_1}_1...y^{\alpha_n}_ndy_1...dy_n=\frac{\Gamma(\alpha_1+1)...\Gamma(\alpha_n+1)}{\Gamma(\alpha_1+...+\alpha_n+n+1)},
\end{equation}
где $\alpha_1,...,\alpha_n\in\mathbb{R}$. Переходя к новым координатам $y_1=x_1^2,...,y_n=x_n^2$, перепишем формулу \eqref{GammaS} в виде
$$
\int\limits_{B^+_1(n)}x_1^{2\alpha_1+1}...x_n^{2\alpha_n+1}dx_1...dx_n=\frac{\Gamma(\alpha_1+1)...\Gamma(\alpha_n+1)}{2^n\Gamma(\alpha_1+...+\alpha_n+n+1)}
$$
или, положив $2\alpha_i+1=\gamma_i$, $i=1,...,n$,
$$
\int\limits_{B^+_1(n)}x_\gamma dx=\frac{\Gamma\left(\frac{\gamma_1+1}{2}\right)...\Gamma\left(\frac{\gamma_n+1}{2}\right)}{2^n\Gamma\left(\frac{n+|\gamma|}{2}+1\right)}.
$$
Тогда, из \eqref{ArBall} и \eqref{Ball1} при $u=1$, используя формулу $\Gamma(z+1)=z\Gamma(z)$ и положив затем $r=1$, получим искомый интеграл
$$
\int\limits_{S^+_1(n)}x^\gamma dS=\frac{\Gamma\left(\frac{\gamma_1+1}{2}\right)...\Gamma\left(\frac{\gamma_n+1}{2}\right)}{2^{n-1}\Gamma\left(\frac{n+|\gamma|}{2}\right)}.
$$
Полученный интеграл обозначается $|S_1^+(n)|_\gamma$:
\begin{equation}\label{PlSh}
|S_1^+(n)|_\gamma
=\int\limits_{S^+_1(n)} x^\gamma dS=
\frac{\prod\limits_{i=1}^n{\Gamma\left(\frac{\gamma_i{+}1}{2}\right)}}{2^{n-1}\Gamma\left(\frac{n{+}|\gamma|}{2}\right)},
\end{equation}
и его также можно найти по формуле 1.2.5 из \cite{LyahovKniga1}, стр. 20, где надо положить $N{=}n$.

\subsection{Многомерный обобщенный сдвиг, весовое сферическое среднее и итерированное весовое сферическое среднее}

В этом пункте мы рассмотрим оператор обобщенного сдвига, действующий по каждой из переменных в части евклидова пространства многих переменных, который будем называть многомерным обобщенным сдвигом и сферическое среднее, где вместо обычного сдвига применяется обобщенный --- весовое сферическое среднее.
Применяя последовательно два раза весовое сферическое среднее, получим итерированное весовое сферическое среднее, которое используется на промежуточных этапах решения задачи о восстановлении функции по ее весовым сферическим средним.

 В  $\mathbb{R}^n_+$ рассматривается многомерный обобщенный сдвиг, отвечающий мультииндексу
$\gamma$ вида
$$^\gamma T^y=\,^{\gamma_1} T_{x_1}^{y_1}...^{\gamma_n}T_{x_n}^{y_n},
$$
где каждый из одномерных обобщенных сдвигов определен выражением
$$
^{\gamma_i} T_{x_i}^{y_i}f(x){=}
\frac{\Gamma\left(\frac{\gamma_i+1}{2}\right)}{\Gamma\left(\frac{\gamma_i}{2}\right)\Gamma\left(\frac{1}{2}\right)}
\int\limits_0^\pi f(x_1,...,x_{i-1},\sqrt{x_i^2+y_i^2-2x_iy_i\cos\alpha_i},x_{i+1},...,x_n)\,\,
\sin^{\gamma_i-1}\alpha_i\,d\alpha_i\,.
$$
На основе многомерного обобщенного сдвига $^\gamma T^y$ конструируется
\textbf{весовое сферическое среднее} функции $f,$ порождённое многомерным обобщённым сдвигом
\begin{equation}\label{05}
M_{f}^{\gamma } (x;t)=M^\gamma_t[f(x)]=\frac{1}{|S_1^+(n)|_\gamma}\int\limits_{S^+_1(n)}\,^\gamma T_x^{t\theta}f(x)
\theta^\gamma dS,
\end{equation}
где $\theta^\gamma{=}\prod\limits_{i=1}^{n}\theta_i^{\gamma_i},$ $S^+_1(n){=}\{\theta{:}|\theta|{=}1,\theta{\in}\mathbb{R}^n_+\}$ --- часть сферы в $\mathbb{R}^n_+$, а
$|S^+_1(n)|_\gamma=\frac{\prod\limits_{i=1}^n{\Gamma\left(\frac{\gamma_i{+}1}{2}\right)}}{2^{n-1}\Gamma\left(\frac{n{+}|\gamma|}{2}\right)}.
$

Весовые сферические средние без обобщенного сдвига изучались в \cite{KiprZas}, \cite{Muravnik0}, \cite{Muravnik1}, стр.117.  В \cite{KiprZas} были получены его свойства и изучен принцип Гюйгенса для волнового уравнения с многими особенностями.  В \cite{Muravnik0}, \cite{Muravnik1} рассмотрены предельные свойства такого среднего. Весовые сферические средние с обычным сдвигом, в случае когда в качестве веса выступал однородный гармонический полином, введены и рассмотрены в \cite{Volchkov1}, стр. 400. Изучению весовых сферические средние, порождённых обобщенным сдвигом, и их применению к решению задач для сингулярных дифференциальных уравнений посвящено большое количество работ, среди которых \cite{LPSh1}--\cite{ShishGert}. Конструкции многомерного обобщенного сдвига и весового сферического среднего представляют собой операторы преобразования (см. \cite{S42} и \cite{S46}).

Отметим простейшие свойства весового сферического среднего (см. \cite{LPSh1}).
\begin{enumerate}
  \item Линейность и однородность: $$M^\gamma_r[af(x)+bg(x)]=aM^\gamma_r[f(x)]+bM^\gamma_r[g(x)],\qquad a,b\in\mathbb{R}.$$
  \item Положительность: если $f(x)\geq 0$, то $M^\gamma_r[f(x)]\geq 0$.
  \item $M^\gamma_r[1]=1$.
%  \item Если $f(x)$ дважды непрерывно дифференцируемая функция, четная по каждой из переменных $x_i$, $i=1,{...},n$, то
% $$
%    (\triangle_{\gamma})_x M_f^\gamma(x;r)=M^\gamma_r(\triangle_{\gamma} f(x)).
 % $$
  \item Пусть ${\bf j}_\gamma(x,\xi){=}\prod\limits_{i=1}^nj_{\frac{\gamma_i-1}{2}}(x_i\,\xi_i),$ где $j_\nu$ определена формулой
  \eqref{FBess}. Справедливо равенство $$M^\gamma_r[{\bf j}_\gamma(x,\xi)]
={\bf j}_\gamma(x,\xi)\,\,j_{\frac{n+|\gamma|-2}{2}}(r|\xi|).$$
    \end{enumerate}

Итерированное весовое сферическое среднее имеет вид (см. \cite{SitShishSemi}, \cite{ShishGert}):
$$
I^{\,\gamma}_f(x;\lambda,\mu)=I^{\,\gamma}_{\lambda,\mu}f(x)=M^\gamma_\lambda M^\gamma_\mu[ f(x)]=$$
$$
=\frac{1}{|S_1^+(n)|_\gamma^2}\int\limits_{S_1^+(n)}\int\limits_{S_1^+(n)}T_x^{\lambda\zeta}T_x^{\mu\xi}[f(x)]\zeta^\gamma
\xi^\gamma dS(\xi) dS(\zeta).
$$

Используя перестановочное свойство 7 обобщённого сдвига, получим, что итерированное весовое сферическое среднее симметрично относительно $\lambda$ и $\mu$:
$$
I^{\,\gamma}_f(x;\lambda,\mu)=I^{\,\gamma}_f(x;\mu,\lambda).
$$
Кроме того, для  $I^{\,\gamma}_f$ имеем очевидные равенства
$$
I^{\,\gamma}_f(x;\lambda,0)=I^{\,\gamma}_f(x;0,\lambda)=M^\gamma_\lambda [f(x)]
$$
и
$$
I^{\,\gamma}_f(x;0,0)=f(x).
$$

В \cite{SitShishSemi} показано, что имеет место равенство, выражающее итерированное сферическое среднее $I^{\,\gamma}_f(x;\lambda,\mu)$ через однократный интеграл от весового сферического среднего $M^\gamma_f(x;r)$. А именно
$$
I^{\,\gamma}_f(x;\lambda,\mu)=\frac{2\,\Gamma\left(\frac{n{+}|\gamma|}{2}\right)}
{\sqrt{\pi}\Gamma\left(\frac{|\gamma|+n-1}{2}\right)}\frac{1}{(2\lambda\mu)^{{n+|\gamma|-2}}}\times
$$
\begin{equation}\label{IterFor}
\times\int\limits_{\lambda-\mu}^{\lambda+\mu}
    \left({(\lambda^2-(r-\mu)^2)((r+\mu)^2-\lambda^2)}\right)^{\frac{n+|\gamma|-3}{2}}\,M^\gamma_f(x;r)\,rdr,
\end{equation}
где $f$ --- интегрируемая с весом $x^\gamma=\prod\limits_{i=1}^nx_i^{\gamma_i}$ функция, чётная по каждой из своих переменных. Если положить $\lambda=\frac{\beta-\alpha}{2}$, $\mu=\frac{\beta+\alpha}{2}$, то равенство \eqref{IterFor} перепишется в виде
$$
I^{\,\gamma}_f\left(x;\frac{\beta-\alpha}{2},\frac{\beta+\alpha}{2}\right)=
$$
\begin{equation}\label{05}
=\frac{\Gamma\left(\frac{n{+}|\gamma|}{2}\right)}{\sqrt{\pi}\Gamma\left(\frac{|\gamma|+n-1}{2}\right)}
\frac{2^{{n+|\gamma|-1}}}{(\beta^2-\alpha^2)^{{n+|\gamma|-2}}}\int\limits_{\alpha}^{\beta}
    \left({(\beta^2-r^2)(r^2-\alpha^2)}\right)^{\frac{n+|\gamma|-3}{2}}M^\gamma_f(x;r)rdr.
\end{equation}
Заметим также, что интеграл в правой части \eqref{IterFor} есть одномерный обобщенный сдвиг (см. формулу \eqref{Sdvig}), таким образом, равенство \eqref{IterFor} при $\lambda>\mu$ примет вид
\begin{equation}\label{Transl}
I^{\,\gamma}_f(x;\lambda,\mu)=\,^\nu{T}_\mu^\lambda M^\gamma_f(x;\mu),\qquad \text{где}\qquad \nu=n+|\gamma|-1.
\end{equation}

В разделе 2 мы рассмотрим дифференциальное уравнение, которому удовлетворяет весовое сферическое среднее --- общее уравнение Эйлера-Пуассона-Дарбу и обобщение теоремы Асгейрссона.

\section{Весовое сферическое среднее, задача Коши и теорема Асгейрссона}

Весовые сферическое средние присутствуют в решениях дифференциальных уравнений в частных производных с оператором Бесселя, используются при решении задач обращения сингулярных интегральных операторов, обобщают сферическое преобразование Радона, имеют приложения в теории выпуклых поверхностей.
Обобщение теоремы И. Стейна (Elias M. Stein) об ограниченности максимальной функции (см. \cite{Stein}), соответствующей весовому сферическому среднему, позволит получать слабые решения общего уравнения Эйлера-Пуассона-Дарбу и доказывать ограниченность интегралов типа свертки с обобщенным сдвигом. Кроме того, используя равенство \eqref{Ball1}, можно рассматривать обобщение максимальной функции Г. Х. Харди (Godfrey Harold Hardy) и Д. И. Литтлвуда (John Edensor Littlewood)(см., например, \cite{Rudin}, стр. 76, \cite{Techl}, стр. 387).

В пункте 2.1 мы рассмотрим задачу Коши с первым ненулевым условием и с нулевым вторым условием для общего уравнения Эйлера-Пуассона-Дарбу. Решение этой задачи выражается через весовые сферические средние (за исключением одного особого случая). Далее, в пункте 2.2 приведено обобщение теоремы Асгейрссона на случай весовых сферических средних и приложение этого обобщения. В последнем пункте 2.3 этого раздела приведены приложения весовых сферических средних.

\subsection{Общее уравнение Эйлера-Пуассона-Дарбу и весовое сферическое среднее}

В этом пункте рассмотрим обобщение классического уравнения Эйлера-Пуассона-Дарбу на случай, когда по всем пространственным переменным действует оператор Бесселя \eqref{Bess}.

Классическое уравнение Эйлера-Пуассона-Дарбу имеет вид
\begin{equation}\label{ClEPD}
\sum\limits^n_{i=1}\frac{{\partial }^2u}{\partial x^2_i}=\frac{{\partial }^2u}{\partial t^2}+\frac{k}{t}\frac{\partial u}{\partial t},\quad u=u(x,t),\quad x\in\mathbb{R}^n,\quad t>0,\quad -\infty<k<\infty.
\end{equation}
%Оператор, действующий по переменной $t$, в \eqref{ClEPD} называется \textbf{оператором Бесселя} и для него принято обозначение (см., например, \cite{kipr}, стр. 3)
%$$
%(B_k)_t=\frac{{\partial }^2}{\partial t^2}+\frac{k}{t}\frac{\partial }{\partial t}.
%$$
Уравнение Эйлера-Пуассона-Дарбу \eqref{ClEPD}  при  $n=1$ впервые было рассмотрено Л. Эйлером (Leonhard Euler) в \cite{Euler} (стр. 227)  и позднее исследовано  в \cite{Poisson} С. Д. Пуассоном (Sim${\rm \acute{e}}$on Denis Poisson), в \cite{Riman} Б. Риманом (Bernhard Riemann)  и в \cite{Darboux} Г. Дарбу (Gaston Darboux) (см. также историю вопроса в \cite{Mizes},  стр. 532). Многомерное уравнения Эйлера-Пуассона-Дарбу \eqref{ClEPD} рассмотрено, например, в \cite{Weinstein3} и \cite{Olevski}. А. Ванштейном (Alexander Weinstein) в статьях \cite{Weinstein0}, \cite{Weinstein1}  решена задача Коши для уравнения \eqref{ClEPD}, где  $k$ в правой части есть произвольное вещественное число и первое условие ненулевое, а второе нулевое. С.А. Терсеновым в \cite{Tersenov} рассмотрена общая задача Коши для \eqref{ClEPD}, где и первое и  второе условия  ненулевые. Кроме того, уравнение Эйлера-Пуассона-Дарбу рассмотрено в \cite{Kuzmin}--\cite{Hairul}.

Пусть \begin{equation}\label{LapBess}
\triangle_\gamma=\sum\limits_{i=1}^{n} (B_{\gamma_i})_{x_i}=\sum\limits_{i=1}^{n} \frac{\partial ^{2} }{\partial x_{i}^{2} } +\frac{\gamma _{i} }{x_{i} } \frac{\partial}{\partial x},\qquad (B_k)_t=\frac{{\partial }^2}{\partial t^2}+\frac{k}{t}\frac{\partial }{\partial t},\qquad k\in\mathbb{R}.
\end{equation}
Весовое сферическое среднее \eqref{05} тесно связано с общим уравнением Эйлера-Пуассона-Дарбу вида
\begin{equation} \label{EPD}
\triangle_\gamma u(x,t)=(B_k)_t\, u , \quad -\infty <k<\infty ,\quad u=u(x,t),\quad x\in \mathbb{R}^n_+ ,\quad t>0,
\end{equation}
Будем рассматривать решения \eqref{EPD}, удовлетворяющие следующим  начальным условиям
\begin{equation}\label{USL}
    u(x,0)=f(x),\qquad \frac{\partial u}{\partial t}\biggr|_{t=0}=0.
\end{equation}
В работе \cite{Fox} (см. также  \cite{CSh}, стр. 243 и \cite{Stellmacher}) и в \cite{LPSh2} рассмотрены различные подходы к решению задачи  Коши  \eqref{EPD}--\eqref{USL} для всех вещественных $k\neq -1,-3,-5,...$. В статье \cite{Barabash}, методом, отличным от методов использующихся в \cite{Fox} и в \cite{CSh} было получено решение \eqref{EPD}--\eqref{USL}  при любых вещественных $k$.

Приведем решения задачи  \eqref{EPD}--\eqref{USL} для различных $k$ (см. \cite{Barabash}). Отметим сначала, что при $k>0$ решение задачи Коши  \eqref{EPD}--\eqref{USL} единственно и не является единственным при $k<0$ (см. \cite{Fox}).

Рассмотрим случай, когда $k=n+|\gamma |-1$. В \cite{LPSh1} доказано, что в этом случае весовое сферическое среднее $M_{f}^{\gamma } (x;t)$ любой дважды непрерывно дифференцируемой функции $f=f(x)$, чётной по каждой из своих независимых переменных $x_{1} ,\ldots ,x_{n} $ удовлетворяет задаче Коши \eqref{EPD}--\eqref{USL}, а именно справедливы равенства
\begin{equation} \label{GrindEQ__6_}
(\Delta _{\gamma } )_{x} M_{f}^{\gamma } (x;t)=(B_{n+|\gamma |-1} )_{r} M_{f}^{\gamma } (x;t){\kern 1pt}   \end{equation}
\begin{equation} \label{GrindEQ__7_} M_{f}^{\gamma } (x;0)=f,\quad \frac{\partial }{\partial t}M_{f}^{\gamma }(x;t)\biggr|_{t=0}=0.
\end{equation}
Формула решения задачи \eqref{EPD}--\eqref{USL} для $k>n+|\gamma |-1$ получается методом спуска и имеет вид
\begin{equation}\label{Sol1}
u^k(x,t)=C(n,\gamma,k)\int\limits_{B_{1}^{+} (n)}[\,^{\gamma}T^{ty}f(x)](1-|y|^2)^{\frac{k-n-|\gamma|-1}{2}}y^\gamma dy,
\end{equation}
где
$$
C(n,\gamma,k)=\frac{\prod\limits_{i=1}^n{\Gamma\left(\frac{\gamma_i{+}1}{2}\right)}\Gamma\left(\frac{k-n-|\gamma|{+}1}{2}\right)}{2^{n}\Gamma\left(\frac{k}{2}\right)}. $$
Переходя в \eqref{Sol1} к сферическим координатам $y=r\theta$, получим  выражение решения $u^k(x,t)$ в виде интеграла от весового сферического среднего $M_f^\gamma$:
\begin{equation}\label{SolSph}
u^k(x,t)=C(n,\gamma,k)|S_1^+(n)|_\gamma\,t^{1-k} \int\limits_0^t M_f^\gamma(x;r)(t^2-r^2)^{\frac{k-n-|\gamma|-1}{2}}r^{n+|\gamma|-1}dr.
\end{equation}
Используя дробный интеграл Эрдейи-Кобера вида
$$
I_{a+;\sigma,\eta}^\alpha f(x)=\frac{\sigma x^{-\sigma(\alpha+\eta)}}{\Gamma(\alpha)}\int\limits_{a}^x(x^\sigma-r^\sigma)^{\alpha-1}r^{\sigma\eta+\sigma-1}f(r)dr,\qquad \alpha>0,
$$
(см. \cite{SKM}, стр. 246, формула 18.1) запишем $u^k(x,t)$ через $I_{a+;\sigma,\eta}^\alpha$:
$$
u^k(x,t)=D(n,\gamma,k)\,I_{0+;2,\frac{n+|\gamma|}{2}-1}^{\frac{k-n-|\gamma|+1}{2}}M_f^\gamma(x;t),
$$
где
$$
D(n,\gamma,k)=\frac{1}{2}C(n,\gamma,k)|S_1^+(n)|_\gamma\,\Gamma\left(\frac{k-n-|\gamma|+1}{2}\right).
$$

Для получения решения при $k{<}n+|\gamma|{-}1$, $k{\neq} {-}1,{-}3,{-}5,{...}$ будем использовать рекуррентные формулы для решения  уравнения \eqref{EPD}. Так, обозначая через $u^{k} =u^{k} (x,t)$ решение уравнения \eqref{EPD} мы имеем две фундаментальные \textbf{рекуррентные формулы}
\begin{equation} \label{GrindEQ__2_} u^{k} =t^{1-k} u^{2-k} , \end{equation}
\begin{equation} \label{GrindEQ__3_} u_{t}^{k} =tu^{k+2} . \end{equation}
Эти рекуррентные формулы позволят при помощи решения $u^{k} $ уравнения \eqref{EPD} получить решение того же уравнения, но с параметром $k+2$ и $2-k$, соответственно. Обе формулы присутствуют в статье Ванштейна \cite{Weinstein0}, но в случае, когда в уравнении \eqref{EPD} по переменным $x_{i} ,\quad i=1,...,n$   действует вторая производная. Формула \eqref{GrindEQ__2_} доказана в \cite{Fox}. Формула \eqref{GrindEQ__3_} доказана в \cite{Barabash}.

Таким образом, в случае когда $k{<}n+|\gamma|{-}1$, $k{\neq} {-}1,{-}3,{-}5,{...}$, используя \eqref{GrindEQ__2_} и \eqref{GrindEQ__3_}, получаем решение \eqref{EPD}--\eqref{USL} (см. \cite{Barabash})
 \begin{equation}\label{34}
u^k=t^{1-k}\,\left(\frac{\partial}{t \partial t}\right)^m(t^{k+2m-1}u^{k+2m}),
\end{equation}
 где $m$ --- положительное целое число  такое, что $k+2m\geq n+|\gamma|-1$,
$$
u^{k+2m}(x,t)=\frac{\prod\limits_{i=1}^n{\Gamma\left(\frac{\gamma_i{+}1}{2}\right)}\Gamma\left(\frac{k+2m-n-|\gamma|{+}1}{2}\right)}{2^{n}\Gamma\left(\frac{k+2m}{2}\right)} \int\limits_{B_{1}^{+} (n)}[\,^{\gamma}T^{ty}f(x)](1-|y|^2)^{\frac{k+2m-n-|\gamma|-1}{2}}y^\gamma dy
$$
   --- решение задачи Коши
\begin{equation}\label{32}
   \Delta_\gamma u^{k+2m}=u^{k+2m}_{tt}+\frac{k+2m}{t}\, u^{k+2m}_t,
\end{equation}
 \begin{equation}\label{33}
u^{k+2m}\left(x,0\right)=\frac{f(x)}{(k+1)(k+3){...}(k+2m-1)},\qquad u^{k+2m}_t\left(x,0\right)=0.
\end{equation}
Поскольку для того, чтобы $u^{k+2m}$ было решением \eqref{EPD}--\eqref{USL}  достаточно, чтобы $f$ имело непрерывную вторую производную, то для получения решения $u^k$ в рассмотренном случае достаточно, чтобы $f$ имела не менее $\frac{1}{2}(n-k+3)$ непрерывных производных.
В этом случае решение рассматриваемой задачи имеет тоже $\frac{1}{2}(n-k+3)$ непрерывных производных.

Запишем, наконец, решение задачи Коши \eqref{EPD}--\eqref{USL}  при  $k{=}{-}1,{-}3,{-}5,{...}$.
Для этого нам потребуется определить
В-полигармоническую функцию. А именно, если функция $f(x)=f(x_1,...,x_n)$ вещественных переменных, определенная в некоторой области пространства $\mathbb{R}^n_+$ имеет непрерывные частные производные до $2m$-го порядка включительно и удовлетворяет всюду в рассматриваемой области В-полигармоническому уравнению
\begin{equation}\label{PolG}
\Delta^m_\gamma f=0,
\end{equation}
где $\Delta_\gamma$ --- оператор, определенный в \eqref{LapBess}, то  $f$ называется \textbf{В-полигармонической функцией} порядка $m$.

Итак, в  случае $k{=}{-}1,{-}3,{-}5,{...}$ решение задачи Коши \eqref{EPD}--\eqref{USL}  $u^k$ существует когда $f$ имеет $\frac{1}{2}(n-k+3)$ непрерывных производных. Однако, частные производные $\frac{\partial^{1-k}u^k}{\partial t^{1-k}}$ этого решения стремятся к бесконечности при $t=0$ со скоростью $\log{t}$ хотя $f(x)$ --- В-полигармоническая порядка $\frac{1-k}{2}$.

 Пусть сначала $k=-1$. Предположим, что $u^{-1}_{tt}(x,0)$ существует. Устремляя  $t$ к $0$ в уравнении
$$
\Delta_\gamma u^{-1}=u_{tt}^{-1}-\frac{u_t^{-1}}{t},
$$
получим, что $\Delta_\gamma u^{-1}(x,0)=0$, то есть $f$ должна быть B-гармонической и решение \eqref{EPD}--\eqref{USL} при $k=-1$ есть $u^{-1}(x,t)=f(x)$.

Переходя далее к $k=-3,-5,...$, получаем, что  решение задачи Коши \eqref{EPD}--\eqref{USL} при $k=-3,-5,...$ дается формулой
 $$
 u^{k}(x,t)=f(x)+\sum\limits_{h=1}^{-\frac{k+1}{2}}\frac{\Delta^h_\gamma f}{(k+1)...(k+2h-1)}\,\frac{t^{2h}}{2\cdot 4\cdot .... 2h},\qquad k=-3,-5,...
 $$

 Отметим, что фундаментальные решения для гиперболического и ультрагиперболического уравнение с оператором Бесселя, примененным по одной переменной получены в \cite{KiprIvanov}, \cite{KiprIvanov1}, по нескольким переменным (кроме временной) в \cite{KiprZas}. Абстрактному уравнению Эйлера-Пуассона-Дарбу посвящены работы \cite{glush1}--\cite{glush6}. Помимо рассмотренного классического подхода к решению уравнения \eqref{ClEPD} и его обобщениям можно применять метод операторов преобразования (см. \cite{Sta1}--\cite{S67}).

 В следующем пункте мы приведем теорему о весовых сферических средних типа теоремы  Асгейрссона и покажем, что с ее помощью также можно получить решение задачи \eqref{EPD}--\eqref{USL}.

\subsection{В-ультрагиперболическое уравнение и обобщение теоремы Асгейрссона}

В этом пункте рассмотри ультрагиперболическое уравнение, в котором вместо вторых производных по каждой переменной действует оператор Бесселя. Такое уравнение будем называть В-ультрагиперболическим. Для него приведем теорему типа теоремы Асгейрссона.

Классическое ультрагиперболическое уравнение вида
\begin{equation}\label{UG}
\Delta_xu=\Delta_yu,\qquad u=u(x,y),\qquad x\in\mathbb{R}^n,\qquad y\in\mathbb{R}^m
\end{equation}
изучалось многими авторами (см. \cite{Blago1}, \cite{Blago2}, \cite{Kosto}, \cite{Berez}, \cite{Owens}, \cite{Radzikowski}, \cite{JohnUG}, \cite{Walter}).
При $n=1$ или $m=1$ \eqref{UG} представляет собой обычное волновое уравнение, описывающая динамическое развитие многих  процессов классической и квантовой физики. К уравнениям  вида \eqref{UG} при $n=m=2$ приводят, например,
 исследования проблемы Гильберта определения в трехмерном декартовом пространстве всех метрик,  геодезическими которых являются прямые (см. \cite{Owens2}); обратная задача дифракции при изучении неоднородности распределения зерен поликристаллических материалов
 (см. \cite{Savelova1}--\cite{Savelova2}); гиперсферическое рентгеновское преобразование, а именно, функции плотности кристаллографических полюсов
  удовлетворяют ультрагиперболическому уравнению с оператором Лапласа-Бельтрами \cite{Nikolayev}. Случай когда в \eqref{UG}  $n>2$ и $m>2$ является важным с математической точки зрения благодаря, в первую очередь, теореме Асгейрссона о среднем значении (см. \cite{Asgeirsson}, \cite{kurantgil2} стр. 475, \cite{john} стр. 84, \cite{Helgason} стр. 318, \cite{Hormander} I, стр. 183).

Теорема Асгейрссона для уравнения   \eqref{UG} имеет следующий вид (см. \cite{kurantgil2} стр. 475).

\begin{teo}\label{tA}
\textbf{Классическая теорема Асгейрссона}. Если $x{\in}\mathbb{R}^n$, $y{\in}\mathbb{R}^n$ и функция $u{=}u(x,y){\in}C^2$ является решением уравнения  \eqref{UG}, то среднее значение функции $u$, взятое при постоянных $x$ в пространстве $\mathbb{R}^n$ переменых $y$ по сфере радиуса $r$, равняется среднему значению $u$, взятому при постоянных $y$ в пространстве $\mathbb{R}^n$ переменных $x$ по сфере того же радиуса $r$:
$$
\int\limits_{|\theta|=r}u(x+\theta,y)dS_\theta=\int\limits_{|\omega|=r}u(x,y+\omega)dS_\omega.
$$
 \end{teo}

Приведем далее результаты из \cite{LPSh1}, \cite{LPSh2}, обобщающие теорему Асгейрссона на случай В-ультрагиперболического уравнения \eqref{EQ22}.
Пусть $n=m'+m''$, $m'$ и $m''$ --- натуральные числа; $\gamma'=(\gamma'_1,...,\gamma'_{m'})$,\,  $\gamma'_i>0$, $i{=}1,...,m'$,
$\gamma''=(\gamma''_1,...,\gamma''_{m''})$,\, $\gamma''_j>0$,\,$j=1,...,m''$\,;  $x\in \mathbb{R}^{m'}_+$, $y\in \mathbb{R}^{m''}_+$, $(x,y)\in\mathbb{R}^n_+=\mathbb{R}^{m'}_+\times \mathbb{R}^{m''}_+$.

\textbf{В-ультрагиперболическое уравнение} или сингулярное ультрагиперболическое уравнение имеет вид
\begin{equation}\label{EQ22}
(\Delta_{\gamma'})_xu=(\Delta_{\gamma''})_yu,\quad u=u(x,y),
\end{equation}
где
$$
(\Delta_{\gamma'})_x=\sum\limits_{i=1}^{m'} (B_{\gamma'_i})_{x_i},\qquad (\Delta_{\gamma''})_y=\sum\limits_{j=1}^{m''}  (B_{\gamma''_j})_{y_j}.
$$

Пусть $u(x,y)\in C^2(\mathbb{R}^n_+)$ и является четной функцией по каждой из своих переменных $x_1,...,x_{m'},y_1,...,y_{m''}$.
 Рассмотрим весовые сферические средние \eqref{05},  взятые по частям поверхностей единичных сфер $S_1^+(m')$,
 $S_1^+(m'')$ в $\mathbb{R}^{m'}_+$ и $\mathbb{R}^{m''}_+$  по каждой из групп переменных $x$ и $y$ с центрами в точках $y\in\overline{\mathbb{R}}\,^{m'}_+$
 и  $z\in\overline{\mathbb{R}}\,^{m''}_+$ соответственно. Для этих весовых сферических средних введем обозначения
$$
(M^{\gamma'}_r u)(x,r;y)=M^{\gamma'}_u (x,r;y)= \frac{1}{|S^+_1(m')|_{\gamma'}}\int\limits_{S^+_1(m')}
T_x^{r\xi}u(x,y)\,\,\xi^{\gamma'}\,dS_\xi,
$$
$$
(M^{\gamma''}_s u)(x;y,s)=M^{\gamma''}_u (x;y,s)=\frac{1}{|S^+_1(m'')|_{\gamma''}}
\int\limits_{S^+_1(m'')}T_y^{s\zeta}u(x,y)\,\,\zeta^{\gamma''}\,dS_\zeta.
$$
Таким образом,  $M^{\gamma'}_u (x,r;y)$ представляет собой весовое сферическое среднее  функции $u$ в  $\mathbb{R}^{m'}_+$ при постоянном  $y=(y_1,\ldots,y_{m''})\in \overline{\mathbb{R}}\,^{m''}_+$\,, а $M^{\gamma''}_u (x;y,s)$ --- весовое сферическое среднее  функции $u$
в  $\mathbb{R}^{m''}_+$ при постоянном  $x=(x_1,\ldots,x_{m'})\in \overline{\mathbb{R}}\,^{m'}_+$.

Составим также общее весовое сферическое среднее функции $u$  по переменным $x$ и $y$
  вида
$$ (M^{\gamma'}_r M^{\gamma''}_s u)(x,r;y,s)=U(x,r;y,s)=$$
$$
 =\frac{1}{|S^+_1(m')|_{\gamma'}\,\,|S^+_1(m'')|_{\gamma''}}
\int\limits_{S^+_1(m'')}\theta^{\gamma''}\,dS(\theta)\,\int\limits_{S^+_1(m')}T_{x,y}^{r\xi,s\zeta}u(x,y)\xi^{\gamma'}\zeta^{\gamma''}dS(\xi).
$$
 Очевидно, что
\begin{equation}\label{EQ23}
M^{\gamma'}_u(x,r;y)=U(x,r;y,0);\quad M^{\gamma''}_u (x;y,s)=U(x,0;y,s).
\end{equation}

\begin{teo}\label{t3}  Если
\begin{equation}\label{EQ24}
m'+|\gamma'|=m''+|\gamma''|
\end{equation}  и
функция $u(x',x'')$ удовлетворяет сингулярному ультрагиперболическому уравнению \eqref{EQ22}, то справедливо равенство:
\begin{equation}\label{EQ25}
M^{\gamma'}_r M^{\gamma''}_s u=M^{\gamma'}_s M^{\gamma''}_r u\,.
\end{equation}
 \end{teo}

Как следствие теоремы \ref{t3} получается результат, являющийся обобщением  теоремы Асгейрссона \ref{tA} на В-ультрагиперболическое уравнение (\ref{EQ22}).

\begin{teo}\label{Ta2}  Пусть функция $u=u(x,y)$ дважды непрерывно дифференцируемая, четная по каждой из своих переменных $x_1,\ldots, x_{m'}$, $y_1,\ldots, y_{m''}$, $(n=m'+m'')$, является решением В-ультрагиперболического уравнения (\ref{EQ22}) и пусть выполняется условие (\ref{EQ24}). Тогда весовое сферическое среднее  функции $u(x,y)$, взятое при постоянных $x{\in}\overline{\mathbb{R}}\,^{m'}_+$ в пространстве $\mathbb{R}^{m''}_+$ по сфере радиуса $r$, равняется весовому сферическому среднему  $u(x,y)$, взятому при постоянных $y{\in}\overline{\mathbb{R}}\,^{m''}_+$ в пространстве $\mathbb{R}^{m'}_+$ по сфере того же радиуса $r$:
 \begin{equation}\label{sled}
 (M^{\gamma'}_u)_{x}(x,y,r)=(M^{\gamma''}_u)_{y}(x,y,r).
 \end{equation}
\end{teo}

Справедливо утверждение обратное по отношению к доказанному следствию, которое представляет собой  обратную теорему Асгейрссона для В-ультрагиперболического уравнения.

\begin{teo} Пусть  $u(x,y){\in}
C^2(\mathbb{R}^{m'}_+{\times} \mathbb{R}^{m''}_+)$ и является функцией четной по каждой из своих переменных $x_1,\ldots,x_{m'},y_1,...,y_{m''}$, $n{=}m'{+}m''$ и пусть для всякой точки $(x,y){\in}
\mathbb{R}^{m'}_+{\times} \mathbb{R}^{m''}_+$ и для всяких неотрицательных $r$ и $s$
выполнено  условие (\ref{sled}). Тогда, если выполнено (\ref{EQ24}), то  функция  $u(x,y)$ удовлетворяет В-ультрагиперболическому уравнению (\ref{EQ22}) в $\mathbb{R}^{m'}_+\times \mathbb{R}^{m''}_+$.
\end{teo}

В \cite{Hormander} на стр. 222 приведено уточнение теоремы Асгейрссона вида:
\begin{teo}\label{utA}(\textbf{Уточнение теоремы Асгейрссона}.) Пусть $u(x,y)$ --- непрерывное решение уравнения
$$
\Delta_xu=\Delta_yu
$$
в некоторой окрестности множества
$$
K=\{(\theta,y)\in \mathbb{R}^{m'}\times\mathbb{R}^{m''}:|\theta|+|\omega|\leq r\},
$$
тогда
$$
\int\limits_{|\theta|=r}u(\theta,0)dS_\theta=\int\limits_{|\omega|=r}u(0,\omega)dS_\omega.
$$
В случае когда $n$ нечетно и отлично от $1$, достаточно предполагать, что $u$ является решением вблизи $\partial K$.
\end{teo}

Справедливо обобщение теоремы \ref{utA} на случай В-ультрагиперболического уравнения  (см. \cite{Suz}):

\begin{teo} Пусть $x{\in} \mathbb{{R}}^{m'}_+,y{\in}\mathbb{{R}}^{m''}_+$, $u{=}u(x,y)$ --- четная   по каждой из своих переменных функция, являющаяся непрерывным  в некоторой  окрестности   множества $K{=}\{ \theta{\in} \mathbb{{R}}^{m'}_+,\omega{\in}\mathbb{{R}}^{m''}_+{:}\,|\theta|{+}|\omega|{=}r\}$ решением В-ультрагиперболического уравнения $(\Delta_{\gamma'})_xu=(\Delta_{\gamma''})_yu,$ и пусть выполнено условие $m'+|\gamma'|=m''+|\gamma''|\ge3$. Тогда справедливо соотношение
$$
\frac{1}{|S_1(m')|_\gamma}
\int\limits_{S_1^+(m')} u(r\theta;0) \,\,\prod\limits_{i=1}^{m'} \theta_i^{\gamma_i}\,dS_\theta=\frac{1}{|S_1(m'')|_\gamma}
\int\limits_{S_1^+(m'')} u(0;r\omega) \,\, \prod\limits_{i=1}^{m''} \omega_i^{\nu_i}\,dS_\omega.
$$
\end{teo}

В следующей теореме приведено решение задачи Коши для общего уравнения Эйлера-Пуассона-Дарбу применением теоремы \ref{Ta2}.
\begin{teo}
 Пусть  $f\in C_{ev}(\mathbb{R}^n_+)$ и $\alpha=\frac{n+|\gamma|-1}{2}$.  Тогда
решение задачи Коши
$$
\left(\frac{\partial}{\partial
t}+\frac{k}{t}\frac{\partial}{\partial
t}\right)u(x,t)=\Delta_{\gamma}u(x,t),\qquad 0<k<1,
$$
$$
u(x,0)=f(x),\quad u_t(x,0)=0
$$
 определяется равенством
$$
u(x,t)=\frac{\Gamma\left(\frac{k+1}{2}\right)}{\Gamma\left(\frac{n+|\gamma|}{2}\right) }
\,t^{1-\delta} {\cal D}^\alpha_r \left(r^{n+|\gamma|-2\over2}M^\gamma_{\sqrt r}\,f(x)\right)(t),
$$
где ${\cal D}^\alpha$ --- левосторонняя дробная  производная Римана-Лиувилля порядка $\alpha$ (см. \cite{SKM}, стр. 44).
\end{teo}

\subsection{О приложениях весовых сферических средних}

Сферические средние тесно связаны с различными методами визуализации, которые представляют большой интерес в различных областях современных исследований, имеющих дело с изображеними в некоторых типах томографических экспериментов, в том числе оптоакустической томографии, термоакустической томографии,  радиолокации и эхолокации (см. \cite{Kuchment0}--\cite{Palamodov1} и др.). Особо важной является задача восстановления функции по ее сферическим средним, решение которой и позволяет получить изображение по данным, полученным в результате сканирования объекта исследования. Чаще всего такой подход встречается в медицине, а именно в таких методах диагностики как  компьютерная томография, магнитно-резонансная томография, ультразвуковое исследование, термография, фотоакустическая томография. Однако, эти же самые методы восстановления изображения используются и для разветки местности при помощи радиолокации и эхолокации. Например, задача восстановления функции по ее сферическим средним является актуальной при получении радиолокационных изображений способом радиолокационного синтезирования апертуры, для решения задач оперативного выявления и оценки фактической радиационной обстановки при помощи позитронно-эмиссионной томографии при исследовании местности методами тепловизорной томографии и др. При этом классическое сферическое среднее далеко не всегда адекватно описывает входные данные изображения, получаемого описанными способами. Во многих задачах (например, в задачах эмиссионной томографии) естественным образом  возникает необходимость использовать весовые сферические средние.

Неотъемлемой частью теории дифференциальных уравнений в частных производных гиперболического типа являются сферические средние (см. \cite{Asgeirsson}--\cite{kurant2}, \cite{IP}--\cite{Ratnakumar}  и др.). В случае, когда рассматриваются уравнения гиперболического типа с оператором Бесселя вместо сферических средних возникают рассмотренные нами весовые сферические средние (см. \cite{KiprZas}--\cite{LPSh3}, \cite{Barabash}, \cite{KiprIvanov}, \cite{Rubin3}  и др.). С точки зрения интегральной геометрии сферические средние рассмотрены в \cite{Volchkov1}, \cite{Helgason}, \cite{Volchkov2}--\cite{Thangavelu}. Большой интерес представляет рассмотрение весовых сферических средних с точки зрения интегральной геометрии.

\vskip 1cm

Во второй части обзора будет рассмотрено обращение весового сферического среднего и его связь с потенциалами Рисса, порожденными обобщенным сдвигом.

\end{document}